\documentclass[11pt]{amsart}

\usepackage{amscd, amssymb}
\begin{document}
\author{Hans U. Boden}
\address{Department of Mathematics, Ohio State University,
Mansfield, OH, 44906}
\email{boden@math.ohio-state.edu}
\author{Christopher M. Herald}
\address{Department of Mathematics, University of Nevada,
Reno, Reno, NV 89557 } 
\email{herald@unr.edu}
 
\title{A Connected Sum Formula  for the SU(3) Casson Invariant}
 
\begin{abstract}
{We provide a formula for the SU(3) Casson invariant 
for 3-manifolds given as the  connected
sum of two integral homology
3-spheres.}
\end{abstract}

\maketitle

\newcommand{\const}{\mbox{const}}
\newcommand{\lto}{\longrightarrow}
\newcommand{\al}{\alpha}
\newcommand{\be}{\beta}
\newcommand{\ga}{\gamma}
\newcommand{\ep}{\epsilon}
\renewcommand{\th}{\theta}
\newcommand{\la}{\lambda}
\newcommand{\om}{\omega}
\newcommand{\si}{\sigma}
\newcommand{\Ga}{\Gamma}
\newcommand{\Om}{\Omega}
\newcommand{\ZZ}{{\mathbb Z}}
\newcommand{\RR}{{\mathbb R}}
\newcommand{\CC}{{\mathbb C}}
\newcommand{\CP}{{\mathbb C}{\mathbb P}}
\renewcommand{\AA}{{\mathcal A}}
\newcommand{\BB}{{\mathcal B}}
\newcommand{\GG}{{\mathcal G}}
\newcommand{\HH}{{\mathcal H}}
\newcommand{\JJ}{{\mathcal J}}
\newcommand{\LL}{{\mathcal L}}
\newcommand{\KK}{{\mathcal K}}
\newcommand{\MM}{{\mathcal M}}
\newcommand{\bG}{{{\mathcal G}_0}}
\newcommand{\NN}{{\mathcal N}}
\newcommand{\OO}{{\mathcal O}}
\newcommand{\UU}{{\mathcal U}}
\newcommand{\hK}{\widehat{K}}
\newcommand{\hQ}{{\widehat{Q}}}
\newcommand{\bB}{{\widetilde{\mathcal B}}}
\newcommand{\bK}{{\widetilde{K}}}
\newcommand{\bM}{{\widetilde{\mathcal M}}}
\newcommand{\bC}{{\widetilde{C}}}
\newcommand{\bU}{{\widetilde{U}}}
\newcommand{\bO}{{\widetilde{\mathcal O}}}
\newcommand{\bL}{{\widetilde{\mathcal L}}}
\newcommand{\wQ}{{\widetilde{Q}}}
\newcommand{\hh}{{\mathfrak h}}
\newcommand{\hhp}{{{\mathfrak h}^\perp}}
\newcommand{\im}{\operatorname{im}}
\newcommand{\hol}{\operatorname{{\it hol}}}
\newcommand{\tr}{\operatorname{\it tr}}
\newcommand{\CS}{\operatorname{{\it CS}}}
\newcommand{\hess}{\operatorname{Hess}}
\newcommand{\grad}{\operatorname{Grad}}
\newcommand{\sym}{\operatorname{Sym}} 
\newcommand{\Spec}{\operatorname{Spec}}
\newcommand{\crit}{\operatorname{Crit}}
\newcommand{\ind}{\operatorname{ind}}
\newcommand{\FF}{{\mathcal F}}
\newcommand{\bbl}{{[\![}}
\newcommand{\bbr}{{]\!]}}
\newcommand{\bA}{{\bbl A \bbr}}
 
\newtheorem{defn}{Definition}
\newtheorem{lem}[defn]{Lemma}
\newtheorem{thm}[defn]{Theorem}
\newtheorem{prop}[defn]{Proposition}
\newtheorem{cor}[defn]{Corollary}
\section{Introduction}
In \cite{bh}, we introduced an invariant $\la_{SU(3)}$ of integral 
homology 3-spheres $X$ defined by appropriately counting  the
conjugacy classes of representations $\varrho:\pi_1 X \to SU(3)$. 
Our main result here is the following theorem.

\begin{thm} \label{connect}  
If $X_1$ and $X_2$ are integral homology 3-spheres, then
\begin{equation} \label{c form} 
\la_{SU(3)}(X_1 \# X_2) = \la_{SU(3)}(X_1) + \la_{SU(3)}(X_2) 
+ 4 \, \la_{SU(2)}(X_1) \; \la_{SU(2)}(X_2),
\end{equation}
where $\la_{SU(2)}$ is Casson's original invariant, 
normalized as in \cite{walker}.
\end{thm}

Even though $\la_{SU(3)}$ is not additive under the connected sum operation, 
the theorem has the following consequence.

\begin{cor} The difference $\la_{SU(3)} - 2\la^2_{SU(2)}$
defines an  invariant of integral homology spheres which is additive
under connected sum. \end{cor}

The proof of Theorem \ref{connect} requires an understanding of 
how certain nondegenerate critical submanifolds of the
(perturbed) Chern-Simons functional contribute to
$\la_{SU(3)}.$ The relevant results here are
Propositions \ref{irred comp} and \ref{red comp}, 
which hold in rather general circumstances.
Before delving into the details, we give a brief introduction 
to 3-manifold $SU(3)$ gauge theory and review the results of \cite{bh}.

Suppose $X$ is a closed, oriented, $\ZZ$-homology 3-sphere
and set $P = X \times SU(3)$.
Denote by $\th$ the trivial (product) connection and by
$d$ the associated covariant derivative.  
Let $$\AA=\{ d+A \mid A \in \Om^1(X;su(3))\}$$ 
be the space of smooth connections in $P$.
The gauge group  $\GG$ of smooth bundle automorphisms $g : P \to P$
acts on $\AA$ by $g \cdot A = g  A g^{-1} + g d g^{-1}$
with quotient $\BB= \AA / \GG,$ the space of gauge orbits of
$SU(3)$ connections.
For the most part we work with the Sobolev completions of
$\AA$ and $\GG$ in the $L^2_1$ and $L^2_2$
norms, respectively, though occasionally we use
the $L^2$ metric on $\AA$.

Denote by  $\Ga_A = \{ g \in \GG \mid g \cdot A = A \}$
the  stabilizer of $A$ in $\GG.$  
The curvature $F_A \in \Om^2(X; su(3))$ is defined for $A \in \AA$
by the formula 
$$F_A  = d A + A \wedge A$$  
and the moduli space $\MM \subset \, \BB$ of flat connections by 
$$\MM = \{ A \in \AA \mid F_A = 0 \} / \GG.$$
For $[A] \in \MM$, $\Ga_A$ is isomorphic to $\ZZ_3$, $U(1)$
or $SU(3)$ because $X$ is a $\ZZ$-homology sphere.
Set $\MM^* = \{ [A] \in \MM \mid \Ga_A = \ZZ_3 \}$ 
and $\MM^r = \{ [A] \in \MM \mid \Ga_A = U(1) \}.$ 
Here, as in \cite{bh}, we call $A \in \AA$ {\it reducible}
if $\Ga_A \cong U(1)$.
Then $\MM$ is the disjoint union $\MM^* \cup \MM^r \cup \{[\th]\}$.

One can also view $\MM$ as the quotient by $\GG$ of the 
critical set of the Chern-Simons functional
\begin{eqnarray*}
&\CS:& \AA \lto \RR \\
&&A \mapsto \tfrac{1}{8 \pi^2} 
\int_X \tr(A \wedge dA + \tfrac{2}{3} A \wedge A \wedge A).
\end{eqnarray*}
First fix a Riemannian metric on $X$ and let 
$* : \Om^p(X;su(3)) \to \Om^{3-p}(X;su(3))$
be the resulting Hodge star operator.
Then define an  inner product on $\AA$ by setting 
$\langle a, b \rangle_{L^{2}} = -\int_{X}\tr(a\wedge *b)$.
Now take the gradient of $\CS$ with respect to the $L^2$ metric on $\AA$ to see that
$$ \grad_A \CS  = -\tfrac{1}{4 \pi^2} * F_A.$$

Consider the self-adjoint elliptic operator $K_A$ which sends 
$(\xi,a) \in \Om^{0 } (X;su(3)) \oplus \Om^1(X;su(3))$ to 
$K_A(\xi,a) = (d_A^* a, d_A \xi - * d_A a).$
Assume   $A$ is flat.  	Then 
$\ker K_A = \HH^0_A(X;su(3)) \oplus \HH^1_A(X;su(3)),$
the space of $d_A$-harmonic $su(3)$-valued (0+1)-forms. 
Choose a path $A_t \in \AA$ with $A_0 = \th$ and $A_1 = A$ and define
${\rm SF}(\th,A)$ to be the spectral flow of the path of 
self-adjoint operators $K_{A_t}.$
If $A$ is reducible, then we can choose the path
so that each $A_t$ has $\Ga_{A_t} = U(1)$ for $t \in (0,1].$
Adjusting by a path of gauge transformations, we can assume that,
for $t \in [0,1], \, A_t \in \AA_{S(U(2) \times U(1))}$, 
the space of connections on $X \times S(U(2) \times U(1)).$ 
Setting $\hh = s(u(2) \times u(1))$ to be the Lie subalgebra of $su(3)$,
it follows that, for $A$ reducible, the spectral flow decomposes as
${\rm SF}(\th,A) = {\rm SF}_\hh(\th,A) + {\rm SF}_\hhp(\th,A)$
according to the splitting $su(3) = \hh \oplus \hhp$,
where $\hhp \cong \CC^2$.

Now $\MM$ is compact and has expected dimension zero 
(since $K_A$ is self-adjoint), 
but it typically contains components of large dimension. 
So that we can work with a discrete space, 
we perturb the Chern-Simons functional using {\it admissible} functions.
These are thoroughly described in Section 2 of \cite{bh}.
Roughly, one alters $\CS:\AA \to \RR$ by adding a gauge-invariant 
function
$h : \AA \to \RR$  of the form $h = \tau \circ \hol_\ell$,
where $\tau : SU(3) \to \RR$ is an invariant function
(usually just the real or imaginary part of trace)
and $\hol_\ell : \AA \to SU(3)$ is the holonomy around 
some loop $\ell \subset X$.
In general circumstances, one must consider sums 
$h=\tau_1 \circ \hol_{\ell_1} + \cdots + \tau_n \circ \hol_{\ell_n}$
where  $\ell_1, \ldots, \ell_n$ are loops in $X$
and $\tau_1,\ldots, \tau_n$ are invariant functions (for analytical 
reasons, one averages these functions over tubular neighborhoods of the 
curves, see \cite{bh} for details). 
Denoting the space of admissible
perturbation functions with respect to this choice
of loops $\ell_1, \ldots, \ell_n$ by $\FF$, by
Definition 2.1 of \cite{bh}, $\FF \cong C^3(\CC,\RR)^{\times n}.$
Each $h \in \FF$ induces a function, also denoted $h$, on $\BB$.

A connection is 
called {\it $h$-perturbed flat} if it is a critical point of $\CS + h.$
Setting $\zeta_h(A) = *F_A - 4 \pi^2  \grad_A h$, 
the moduli space of $h$-perturbed flat connections is defined to be
$$\MM_h = \zeta_h^{-1}(0)/\GG.$$ 
We denote by $\MM^*_h$ (and $\MM^r_h$) the subset of gauge orbits of
irreducible (reducible, respectively) perturbed flat connections.
 
Perturbing only changes the flatness equation in a small neighborhood
of the supporting loops $\ell_i$. For example, when 
$h= \tau \circ \hol_\ell$,
every perturbed flat connection $A$ is actually flat outside
a small tubular neighborhood of $\ell$.
In  general if $h=\sum_{i=1}^n \tau_i \hol_{\ell_i}$,
then the same is true outside the union of small tubular neighborhoods of each $\ell_i.$
We showed in Section 3 of \cite{bh}
that there exist loops $\ell_1, \ldots, \ell_n$ in $X$ such that, for 
generic  small $h \in \FF$, $\MM^{*}_{h}$ and 
$\MM^{r}_{h}$ are compact 0-dimensional 
submanifolds of $\BB^{*}$  and $\BB^{r}$ consisting of gauge orbits 
that satisfy a cohomological regularity condition  
(see Definition \ref{regular} below and Theorem 3.13 of \cite{bh}).
Moreover, if $A$ is
$h$-perturbed flat, then there is a flat connection
$\widehat{A}$ near  $A$ 
(cf. Proposition 3.7, \cite{bh}).

\begin{prop}  \label{mainer}
For generic $h$ sufficiently small, 
the quantity
$$\la_{SU(3)}(X) : = \sum_{[A] \in \MM^*_h} (-1)^{{\rm SF}(\th, A)}
- \frac{1}{2}\sum_{[A] \in \MM^{r}_h}(-1)^{{\rm SF}(\th, A)}
({\rm SF}_{\hhp}(\th, A) - 4 \CS(\widehat{A}) + 2)$$
defines an invariant of integral homology 3-spheres $X$
called the {\bf Casson SU(3) invariant}.
\end{prop}

In reference to the second sum, only the difference  
${\rm SF}_\hhp(\th,A) - 4 \CS(\widehat{A})$ 
is well-defined on the gauge orbit $[A]$; 
each term individually depends on the choice 
of representative for $[A].$
It is proved in \cite{bh} that  the above formula for $\la_{SU(3)}(X)$
is independent of the choice of $h,$ Riemannian metric, 
and orientation of $X.$

\section{The Gluing Construction and Point Components}
Theorem \ref{connect} is proved by gluing together 
perturbed flat connections on $X_1$ and $X_2$.
For $i=1,2,$ set $P_i = X_i \times SU(3)$ and denote by $\th_i$ the 
trivial 
connection in $P_i$. 
Choose $h_i$ a generic sufficiently small admissible perturbation function 
so that $\MM_{h_i}(X_i)$, the moduli space of perturbed flat  
connections in $P_i$, is {\it regular}
according to the following definition. 

We first introduce some notation.
Given a smooth function $h:\AA\to \RR$,
the Hessian of $h$ at $A$ is the map $\hess_A h: \Om^1(X;su(3)) \to \Om^1(X;su(3))$  
defined in terms of the $L^2$ metric  by  
$$\langle \hess_A f(a), b \rangle_{L^2} = 
\left.\frac{\partial^2}{\partial s \partial t} h(A + sa + tb) \right|_{s,t=0}.$$

\begin{defn}  \label{regular}
Suppose $X$ is a $\ZZ$-homology 3-sphere,
$P = X \times SU(3)$,
$h : \AA \to \RR$ is an admissible perturbation function and $A$
is an $h$-perturbed flat connection. Introduce the
operator $*d_{A,h} = *d_A - 4 \pi^2 \hess_A h $ on $\Om^{1}(X; su(3))$
and define the {\bf deformation complex} to be
\begin{eqnarray}
\Om^{0}(X;su(3)) \stackrel{d_{A}}{\lto}
\Om^{1}(X;su(3)) \stackrel{*d_{A,h}}{\hspace{0mm}\lto}
\Om^{1}(X;su(3)) \stackrel{d_{A}^{*}}{\lto}
\Om^{0}(X;su(3)).
\end{eqnarray}
Define groups
 $H^0_A(X;su(3)) = \ker d_A$ (the   Lie algebra of the stabilizer subgroup $\Ga_A$)
and $H^{1}_{A,h}(X;su(3)) = \ker *d_{A,h} / \im d_A$.
A point $[A] \in \MM_h$ is called {\bf regular} if 
$H^1_{A,h}(X, su(3)) = 0,$ 
and a subset $S \subseteq \MM_{h}$ is regular if this condition holds 
for all $[A] \in S.$ 
\end{defn}

The procedure outlined in \S 7.2.1 of \cite{donaldson kronheimer}
constructs a nearly anti-self-dual connection on $X_1 \# X_2$
given anti-self-dual connections
$A_1$ and $A_2$ on 4-manifolds $X_1$ and $X_2$.
A key step is to approximate $A_i$ by a connection that is
flat in a small neighborhood of the basepoint $x_i \in X_i$.
We use a similar (but simpler) procedure to construct perturbed
flat connections on the connected sum of two 3-manifolds.
We first review the construction for $X_1 \# X_2$, 
then construct the bundle $P_1 \# P_2$ and connection
(see also \cite{wli}).

Given basepoints $x_i \in X_i$ and small, 3-balls $B_i$ containing $x_i,$
set  $ \dot{B}_i = B_i \setminus \{x_i \}$
and $\dot{X}_i = X_i \setminus \{x_i \}.$  We take the metric to be 
flat on $B_i$.  
Choose an orientation reversing isometry 
$f: \dot{B}_1 \to \dot{B}_2$
of the deleted neighborhoods and define
$X_1 \# X_2 = \dot{X}_1   \cup \dot{X}_2/\sim,$
where $x \sim f(x) $ for $x \in \dot{B}_1$.

Now suppose 
$h_1 = \sum_{j=1}^{n_1}\tau_{1,j} \hol_{\ell_{1,j}}$ 
and $h_2 = \sum_{j=1}^{n_2}\tau_{2,j} \hol_{\ell_{2,j}}$ are admissible
perturbations on $X_1$ and $X_2$, respectively.
We can choose $x_i$ and $B_i$ 
so that $\ell_{i,j}$ misses $B_i$ for all 
$j=1, \ldots n_i$ and each $i=1,2$.
Thus, if $A_i$ is an $h_i$-perturbed flat connection on $X_i$,
its restriction to $B_i$ is flat and parallel translation
by $A_i$ defines a trivialization of $P_i|_{B_i}$ in which the
connection is also trivial.

Using these trivializations, we can extend
any isomorphism $\si : (P_1)_{x_1} \to (P_2)_{x_2}$ to
an isomorphism of $P_1|_{{B}_1} \to P_2|_{{B}_2}$.
We then construct the bundle $P_1 \# P_2$ 
by gluing  $P_1$ and $P_2$
by identifying  $P_1|_{\dot{B}_1}$ and 
$P_2|_{\dot{B}_2}$. Since the restriction of $A_i$ to
$B_i$ is trivial, we can also glue $A_1$ and $A_2$
to obtain the connection
$A_1 \#_\si A_2$ on $ P_1 \# P_2$.
Of course, $P_1 \# P_2 \cong X \times SU(3)$ is independent of
$\si$ even though $A_1 \#_\si A_2$ is not, in general.

Since the loops
 $\ell_{i,j}$ do not intersect the balls $B_i$, setting $h_0=h_1+h_2$
 defines an admissible perturbation on $X=X_1 \# X_2$.
If $A$ is an $h_0$-perturbed flat connection  on $X,$ 
then restricting $A$ to each side of the connected sum, 
shows that 
$A$ is gauge equivalent to one the form $A_1 \#_\si A_2$ for some $A_1, A_2$ 
and $\si$ as above.
Moreover, $A_1 \#_\si A_2 $ and $A_1 \#_{\si'}A_2$ 
are gauge equivalent
if and only if $\si$ and $\si'$ are in the same 
$\Ga_{A_1} \times \Ga_{A_2}$ orbit in $SU(3).$

Observe that $\MM_{h_0}(X)$ is not regular,
even though both $\MM_{h_1}(X_1)$ and $\MM_{h_2}(X_2)$ are.
In fact, the gauge orbit $[A_1 \#_\si A_2]$ is  isolated in 
$\MM_{h_0}(X)$
if and only if $A_i = \th_i$ for $i=1$ or $2$. 
In that case, $[A_1 \#_\si A_2]$ is independent of $\si$
and so we drop the subscript and simply write
$[A_1 \# \th_2]$ or $[\th_1 \# A_2].$ 
 
Since $\MM_{h_0}(X)$ is not regular, 
one cannot compute
$\la_{SU(3)}(X)$ 
from Proposition \ref{mainer} without further
perturbing the flatness equations.
A method for doing this is presented in the next section, but first 
we explain the special role played by
connections of the form $A_1 \# \th_2$ and $\th_1 \# A_2$.
By a Mayer-Vietoris argument,
the gauge orbits $[A_1 \# \th_2]$ and
$[\th_1 \# A_2]$ in $\MM_{h_0}(X)$ are regular whenever
$[A_1] \in \MM_{h_1}(X_1)$ and $[A_2] \in \MM_{h_2}(X_2)$ are regular.
If $C \subset \MM_{h_0}(X)$ is a point component,
then   either $C = \{[A_1 \# \th_2] \}$ or   $C = \{ [\th_1 \# A_2] \}.$

It is well-known that for irreducible connections,
the spectral flow is additive with respect to connected sum.
Specifically,
if $\th=\th_1 \# \th_2$ and $A=A_1 \#_\si A_2$ where $A_1$ and 
$A_2$ are irreducible connections on $X_1$ and $X_2$, respectively, then
\begin{equation} \label{splitit}
{\rm SF}_{X}(\th,A) = {\rm SF}_{X_1}(\th_1,A_1) + {\rm SF}_{X_2}(\th_2,A_2). 
\end{equation}
(For proofs of this statement and the next in the $SU(2)$ setting, 
see Lemmas 2.2.1 and 2.2.2 in \cite{wli}.)
The next result
treats the case
when $A_1$ or $A_2$ is trivial and determines
the contribution of point components to $\la_{SU(3)}(X_1 \# X_2).$

\begin{lem}  \label{addit}  Set $\th = \th_1 \# \th_2$ 
and suppose that $A_i$ is a nontrivial, $h_i$-perturbed flat
$SU(3)$ connection on $X_i$ for $i=1,2$.
In parts (ii) and (iii), assume further that $A_i$
is reducible and that $\widehat{A}_i$ is the reducible flat connection on 
$X_i$
close to $A_i$ for $i=1,2.$ Then
\begin{enumerate}
\item[(i)] 
${\rm SF}_X(\th, A_1 \# \th_2)  = {\rm SF}_{X_1}(\th_1, A_1)\;$ and
$\; {\rm SF}_X(\th, \th_1 \# A_2)  = {\rm SF}_{X_2}(\th_2, A_2).$ 
\item[(ii)] 
${\rm SF}_{X, \hhp}(\th, A_1 \# \th_2)  = {\rm SF}_{X_1, \hhp}(\th_1, A_1)\;$ and 
$\; {\rm SF}_{X, \hhp}(\th, \th_1 \# A_2)  = {\rm SF}_{X_2, \hhp}(\th_2, A_2).$
\item[(iii)]

$\CS_X(\widehat{A}_1 \# \th_2) = \CS_{X_1}(\widehat{A}_1) \;$ and 
$\; \CS_X(\th_1 \# \widehat{A}_2) = \CS_{X_2}(\widehat{A}_2 ).$
\end{enumerate}
\end{lem}
 
Using Lemma \ref{addit} and summing over the set 
$$\MM^0_{h_0}(X) = \{ [A] \in \MM_{h_0}(X) \mid A=A_1 \# \th_2 \hbox{ or 
} 
A = \th_1 \# A_2  \}$$
of point components of $\MM_{h_0}(X),$ we see that

\begin{align*}  
{\sum_{[A] \in \MM^{0,*}_{h_0}(X)}} \!\!\!\!\!\!
& (-1)^{{\rm SF}(\th, A)}
- \tfrac{1}{2} \!\!\!\!\!\!
{\sum_{[A] \in \MM^{0,r}_{h_0}(X)}} \!\!\!\!\!\!
(-1)^{{\rm SF}(\th, A)}
({\rm SF}_{\hhp}(\th, A) - 4 \CS(\widehat{A}) + 2) \\ 
&= \!\!\!\! \!\! \sum_{[A_1] \in \MM^*_{h_1}(X_1)} 
\!\!\!\! \!\! (-1)^{{\rm SF}(\th_1, A_1)} - \tfrac{1}{2} \!\!\!\!\!\! 
 \sum_{[A_1] \in \MM^r_{h_1}(X_1)} \!\!\!\!\!\!
 (-1)^{{\rm SF}(\th_1, A_1)} ({\rm SF}_{\hhp}(\th_1, A_1) 
 - 4 \CS(\widehat{A}_1) + 2)  \\  
&\qquad +  \!\!\!\! \!\! \sum_{[A_2] \in \MM^*_{h_2}(X_2)}  
\!\!\!\!\!\! (-1)^{{\rm SF}(\th_2, A_2)} - \tfrac{1}{2} \!\!\!\! \!\!
 \sum_{[A_2] \in \MM^r_{h_2}(X_2)} \!\!\!\!\!\!
(-1)^{{\rm SF}(\th_2, A_2)} ({\rm SF}_{\hhp}(\th_2, A_2) 
- 4 \CS(\widehat{A}_2) + 2)&\\ 
&= \la_{SU(3)}(X_1) + \la_{SU(3)}(X_2).\\
\end{align*} 
Thus, the point components in $\MM_{h_0}(X)$
give rise to the first two terms on the right hand side of
formula (\ref{c form}).

\section{Higher Dimensional Components}
In this section, we study connected components $C$
of $\MM_{h_0}(X)$ with $\dim C >0$ and   
analyze their contribution to $\la_{SU(3)}(X_1 \# X_2).$
Here and elsewhere in this section,
$h_0=h_1 + h_2$ is the perturbation from the previous section
obtained by perturbing over $X_1$ and $X_2$ separately.
Suppose $C$ is such a component
and suppose $[A_1 \#_\si A_2] \in C$.
Then, since $\MM_{h_1}(X_1)$ and $\MM_{h_2}(X_2)$
are both regular, we obtain an explicit description
of $C$ as the double coset
space of $SU(3)$ by $\Ga_{A_1}$ and $\Ga_{A_2}$.

We also introduce
the {\it based} gauge group 
$\bG = \{ g \in \GG \mid g_{x_0} = 1 \}$,
where $x_0 \in X$ is a fixed
basepoint. Set $\bB = \BB/\bG$, the space of based gauge
orbits of connections, and $\bM_h = \zeta_h^{-1}(0) / \GG_0$,
the based perturbed flat moduli space.
Using the gluing construction, it is not difficult to see that
$\bM_{h_0}(X_1 \# X_2) = \bM_{h_1}(X_1) \times \bM_{h_2}(X_2)$.

The projection $\pi:\bB \to \BB$ has fiber
modeled on $SU(3)/\Ga_A$ over $[A].$  
The two fiber types relevant here are
$PU(3) = SU(3)/\ZZ_3$
and the homogeneous 7-manifold $N$  obtained as the space of
left cosets of the $U(1)$ subgroup
\begin{equation}\label{U1 subgroup}
 \left\{ \left. \left( 
\begin{array}{ccc}
u&0&0\\
0&u&0\\
0&0&u^{-2} \end{array}
\right)  \right| u \in U(1) \right\}\end{equation}
of $SU(3)$. 
From now on, since we will be dealing almost exclusively
with connections on $X=X_1 \#X_2,$ we write $\MM_{h}$ for $\MM_{h}(X)$. 
The following proposition summarizes what
we now know about the components $C \subset \MM_{h_0}$ with $\dim C >0.$
 
\begin{prop}   
Suppose $C= \{[A_1 \#_\si A_2] \mid \si \in 
\Ga_{A_1} \backslash SU(3)/\Ga_{A_2}\}$
is a connected component of $\MM_{h_0}$,
where both $A_1$ and $A_2$ are nontrivial (so $C$ is not a point 
component).
\begin{enumerate}
\item[(i)]  If  $A_1$ or $A_2$ is irreducible, then 
$C$ is a smooth submanifold of $\BB^*$ with
$ C \cong PU(3)$ if  $A_1$ and $A_2$ are both irreducible, and 
$C \cong N $ if $A_1$ or $A_2$ is reducible. 
\item[(ii)] If  both $A_1$ and $A_2$ are reducible, then  
$\bC \cong N \times N$ 
is a smooth submanifold of $\bB $, where $\bC$ is the preimage of $C$
under the projection $\pi:\bB \to \BB$.
\end{enumerate}
In (i), the component $C$ is nondegenerate, that is, the
Hessian of $\CS + h_0$ is nondegenerate in the normal directions to $C$.
In (ii),  the same is true of $\bC.$  
\label{component types}
\end{prop}

Obviously $h_0 \in \FF,$ and for generic  $h$  near $h_0,$
the moduli space $\MM_{h}$ will be  regular and 
every  $[A] \in \MM_{h}$ will be close to some $[A_0] \in \MM_{h_0}.$ 
Moreover, for components $C$ of type (i), the restriction
$h|_C$ will generically be a Morse function.
To see this, consider the bundle $E$ over $\FF \times C$ obtained
from $TC \to C$ by pullback under $\FF \times C \to C$.
Define a section $s:\FF \times C \to E$
by setting $s(h,[A]) = \grad _{[A]} (h|_C) $.
The abundance condition implies $s$ is 
a submersion, and thus we have an open set $V$ in $\FF$
containing $h_0$ and a subset   $V' \subset V$
of second category such that $h \in V'$ implies $h|_C$ is  Morse. 

Using such $h$, we can evaluate the contribution to 
$\lambda_{SU(3)}(X_1 \# X_2)$ of the critical
points in $\MM_{h}$ arising from each component $C \subset \MM_{h_0}$.
For components of type (i), we apply the following  lemma. 
Although the result is well-known, we include
a proof because we could not find one in the literature.
This proof will later be generalized to establish Lemma \ref{bott morse 
pert}, 
an equivariant version of this result which is new, as far as we know.

\begin{lem}\label{morsepert}  
Suppose $C\subset \MM_{h_0}^{*}$ is a 
nondegenerate critical submanifold
and $f$  is an admissible function with $f|_{C}$ Morse.
Set $h_t = h_0 + t f$ for $t$ small.
Then there is an open set $U \subset \BB^{*}$ containing $C$
and an $\ep>0$ such that, for every $0<t<\ep$, 
$\OO_t :=\MM_{h_t}\cap U$ is a 
regular subset of $\MM_{h_t}$ with a natural bijection
$\varphi_{t}:\crit (f|_{C}) \to \OO_t$.
Given  a smooth family of connections $A_t$ with  
$[A_0] \in \crit(f|_C)$ and
$[A_t] = \varphi_t([A_0])$   for $0<t<\ep,$ then
\begin{equation} \label{sfclaim}
{\rm SF}(\th,A_t) = {\rm SF}(\th,A_0)  + \ind_{[A_0]}(f),
\end{equation}
where $\ind_{[A_0]}(f) $ is the Morse index of 
the critical point $[A_0]$ with respect to the function  $f|_C$.
\end{lem}
\begin{proof}
 
We begin by introducing some notation and recalling some basic material 
from
\cite{taubes} and \cite{bh}. 
Let $\JJ$ be the trivial bundle over $\AA \times \FF$ with fiber 
$\Om^{0+1}(X;su(3))$. Impose the $L^{2}$
pre-Hilbert space structure on the fibers and consider the 
smooth subbundle $\LL\subset \JJ\mid _{\AA^* \times \FF}$
 whose fiber above $(A,h)$ is
$$\LL_{A,h} =\{(\xi, a)\in \JJ_{A,h} \mid  \xi=0, d_{A}^{*}a=0 \}.$$  
The bundle  $\LL$ over   $\AA^* \times \FF$ is $\GG$-equivariant 
and hence descends to give a  bundle, also denoted by $\LL$,  
over $\BB^{*}\times \FF$, which we regard as the tangent bundle to 
$\BB^{*}$  
with the $L^{2}$ metric as opposed to a Sobolev metric.

Recall the operator $K_A$ on $\Om^{0+1}(X;su(3))$ defined by 
$K_A(\xi,a) = (d^*_A a, d_A \xi - *d_A a)$.
It can be extended to give an operator $K:\JJ \to \JJ$ by setting
$$K_{A,h}(\xi, a)=\left(d_{A}^{*}a, d_{A}\xi - *d_{A,h}a\right)
= \left(d_{A}^{*}a, d_{A}\xi - *d_{A} a + 4 \pi^2 \hess_A h (a)\right).$$  
Then $K_{A,h}$ is a closed, essentially self-adjoint Fredholm operator
with dense domain, depending smoothly on $A$ and $h$.  It has discrete
spectrum with no accumulation points, and each eigenvalue has finite
multiplicity. If $A$ is $h$-perturbed flat,
then $K_{A,h}$  respects the splitting $\JJ = \LL' \oplus \LL$ 
where $\LL' = \Om^0 \oplus {\rm Im}(d_A:\Om^0 \to \Om^1)$. 
 
\bigskip \noindent
{\it Remark. \;} 
Note that $K_{A,h}$ as defined here differs from the operator
used in \cite{bh}. However, the formula for $\la_{SU(3)}(X)$
is the same, because changing the sign of $*d_A$ in $K$ is
equivalent to changing the orientation of the 3-manifold,
and it is proved in \cite{bh} that $\la_{SU(3)}(-X) = \la_{SU(3)}(X).$
\bigskip

We now introduce a  closely related operator on $\LL$.
Let $\pi_{A,h}:  \JJ_{A,h}  \to \LL_{A,h}$ be the 
$L^2$-orthogonal projection and let $\hK_{A,h}$ be 
the  operator on $\LL_{A,h}$ obtained by restricting 
$\pi_{A,h} \circ   K_{A,h}$. 
For paths in $\FF \times \BB^*,$ the spectral flow
of $K_{A,h}$ and $\hK_{A,h}$ are identical.

Let 
\begin{equation} \label{lam}
\la_0 = \min \{ |\la| \;  \mid \la \neq 0, \la 
\in \Spec(\hK_{A_0,h_0}) \mbox{ for } [A_0]\in C\}.
\end{equation}
Choose open neighborhoods $U\subset \BB^{*}$ of $C$ and $V \subset \FF$
of $h_0$ small enough so that
$([A],h) \in U \times V$ implies $\la_0/2 \not\in \Spec(\hK_{A,h})$.  
Over $U \times V$ it is possible to decompose $\LL$ into
$\LL_0\oplus \LL_1$ where 
\begin{equation} \label{big and small}
\LL_0 = \bigoplus_{|\la| < \la_0/2} E_\la
\quad \quad \hbox{and} \quad \quad 
\LL_1 = \bigoplus_{|\la| > \la_0/2} E_\la.
\end{equation}
Here $\la \in \Spec (K_{A,h})$ is an eigenvalue and $E_\la$ is its 
eigenspace.

Let $p_i : \LL \to \LL_i$  be
the projection and choose $\ep > 0$ so that $h_t \in V$ for $t \in 
(-\ep,\ep)$.  
For $i=0,1,$ define $$\psi_i: U \times (-\ep,\ep) \to \LL_i$$
by setting $\psi_i([A],t) = p_i(\zeta_{h_t}(A))$. 
(Recall that $\zeta_h(A) = *F_A - 4 \pi^2 \grad_A h $.)
A standard argument shows that $\psi_1$ is a submersion 
along $C\times\{0\}$ and so, by the Inverse Function Theorem,
for $U$ and $\ep$ small enough, $\psi_1^{-1}(0)$ is a  
submanifold of $U \times (-\ep,\ep) $ parameterized by a $C^3$ function 
$\Phi:C\times (-\ep, \ep) \to U \times  (-\ep, \ep)$
of the form $\Phi([A],t) = (\phi_t([A]),t)$, where $\phi_t:  C \to U$ is 
smooth.
 
Consider part of the parameterized moduli space
$W={\displaystyle \bigcup_{t \in (-\ep,\ep)}} \MM_{h_t} \times \{t\}$ 
defined by
$$W_\ep= \{([A],t) \mid [A] \in U, \, -\ep < t< \ep, \, 
\zeta_{h_t}(A)=0\}.$$
Then $W_\ep$ is the image under $\Phi$ of the zero set
of the map $Q$ from $C\times (-\ep, \ep)$ to 
$\LL_0$ defined by $Q=\psi_{0}\circ \Phi$.  
This zero set is not cut out
transversely since $\MM_{h_0}$ is not regular along $C.$
We expand $Q(x,t)$ about $t=0$ for $x \in C$.  
For clarity we are
using $x$ instead of $[A]$ to denote gauge orbits.    
Since $x \in C,$   $\zeta_{h_0}(x)=0$ and  we have
$$\zeta_{h_t}(\phi_t(x)) =  t \hess_{x} (\CS + h_{0})  
\left(\left. \tfrac{ d \phi_t(x)}{dt} \right|_{t=0} \right) 
- 4 \pi^2 t \grad_x  f + O(t^2).$$ 
It then follows that
\begin{eqnarray*}
Q(x,t) &=& 
  p_0 \left( \zeta_{h_t}(\phi_t(x)) \right)\\
&=&   p_{0}\left[ t \hess_{x} (\CS + h_{0})  
\left(\left. \tfrac{ d \phi_t(x)}{dt} \right|_{t=0} \right) 
- 4 \pi^2  t \grad_x f \right]+ O(t^2) \\
&=&  -4 \pi^2 \, t \,  p_0  (\grad_x f)  + O(t^2).
\end{eqnarray*}
This last step follows since $p_0$ is the projection onto the kernel
of the Hessian of $\CS+h_0$.
Thus the function $Q/t$ extends to a $C^2$ function  
$\hQ:C\times (-\ep, \ep) \lto \LL_0 $ defined by
$$\hQ(x,t)=\left\{ \begin{array}{ll}
\ \ 
Q(x,t)/t & \hbox{if } t\neq 0\\
-4 \pi^2 p_0 (\grad_x f) & \hbox{otherwise.}
\end{array} \right.$$
Obviously, for $t \neq 0,$ the zero set of $\hQ$ coincides with that of 
$Q$.
Moreover, the restriction of $\hQ$ to $C\times \{0\}$ is 
transverse to the zero section of $\LL_0$,
since by hypothesis $f|_{C}$ is a Morse function.  Therefore, 
for $\ep$ small enough, $\hQ^{-1}(0)$
is a smooth, 1-dimensional submanifold
of $C \times (-\ep,\ep)$ which  intersects $C\times \{0\}$ transversely 
and 
$$\hQ^{-1}(0) \cap (C \times \{0\}) = \crit(f|_C).$$
Following this product cobordism gives a natural bijection 
$\varphi_t:  \crit(f|_C) \to \OO_t$.

To prove (\ref{sfclaim}), let $[A_0] \in \crit(f|_{C})$ and
denote by $A_{t}$ a differentiable family of connections 
representing the path of orbits $\varphi_{t}([A_0])$.
Consider the differentiable family of closed, essentially self-adjoint 
Fredholm
operators $K(t):= \hK_{A_{t},h_t}$. (Here we could equally well
work with the path $K_{A_{t},h_t}$ of operators on $\JJ$ since we
are only concerned with the behavior of the small eigenvalues.)

The eigenvalues of $K(t)$ of modulus less than $\la_0$ vary
continuously differentiably in $t$, and their derivatives 
at $t=0$ are given by the eigenvalues of 
$\left. \frac{\partial K(t)}{\partial t}\right|_{t=0}$
restricted and projected to $\ker K(0) = \ker \hK_{A_0,h_0}$  
(see Theorem II.5.4  and Section III.6.5 of \cite{kato}).
However, one can see directly that the restriction of
$ p_0 \left( \left. \frac{\partial K(t)}{\partial t}\right|_{t=0}\right)$
to $\LL_0$ agrees with $\hess_{[A]} (f|_{C})$ and this completes the 
proof.
\end{proof}

The next result applies to components of type (i)
and determines their contribution to $\la_{SU(3)}(X).$ 
It is an immediate consequence of Lemma \ref{morsepert}.

\begin{prop}  \label{irred comp}
Suppose $C \subset \MM^*_{h_0}$ is a nondegenerate critical
submanifold and $[A] \in C.$ Then the contribution of $C$ to 
$\la_{SU(3)}(X)$ is $(-1)^{{\rm SF}(\th,A)} \chi(C).$ 
\end{prop}

Next we develop similar results for components $C$ of type (ii).
In this case, since $C$ is not smooth,
we work equivariantly on $\bC$,
which has a natural $SU(3) \cong \GG/\bG$ action.
First, we introduce a relevant definition.

\begin{defn}    \label{equiv morse}
Suppose $G$ is a compact Lie group acting smoothly on
a compact manifold $Y$. Then a smooth $G$-invariant function $f:Y \to 
\RR$ is 
called {\bf equivariantly Morse} if its critical point set 
$\crit(f)$ is a union of  orbits isolated in $Y/G$ and along any 
such orbit the Hessian of $f$ is nondegenerate in the normal directions.
\end{defn}
\noindent
Note that an equivariantly  Morse function is not necessarily Morse,
though it is always Bott-Morse.

Let $\bC^*$  and $\bC^r$  be the preimages  of $C^*$  and $C^r$  
under the projection $\pi:\bB \to \BB$.
They determine  a stratification  $\bC = \bC^* \cup \bC^r$
 given by orbit type.
We denote by $\bA \in \bB$ the $\bG$ orbit
of $A \in \AA$. 
Observe that $\Ga_A \cong \ZZ_3$ for $\bA \in \bC^*$ and
$\Ga_A \cong U(1)$ for $\bA \in \bC^r.$
This latter isomorphism endows $\nu(\bC^r),$
the normal bundle  of $\bC^r$ in $\bC$, with a natural  $U(1)$ action. 
Every $h \in \FF$ defines an invariant function on $\bC$
by restriction. If $h$ is equivariantly  Morse and
$\tau \subset \bC$ is an open, $SU(3)$ invariant tubular neighborhood of 
$\bC^r$, then 
the induced functions  $ (\bC^* \setminus \tau)/SU(3) \to \RR$
and $ \bC^r/SU(3) \to \RR$  obtained by restricting and 
passing to the quotient are both Morse functions
with only finitely many critical points.

We now prove that generic $h \in \FF$ induce
equivariantly  Morse functions on $\bC$.
This is achieved in two steps.	
First, let $\xi$ be the bundle over $\FF \times \bC^r$ obtained
by pulling back  
the bundle $T \bC^r \oplus \sym(\nu)$ 
under $\FF \times \bC^r \to \bC^r$, where  
$\sym(\nu)$ is the bundle of $U(1) $ equivariant
symmetric bilinear forms on $\nu(\bC^r)$.  Define a section 
$s:\FF \times \bC^r \lto \xi$ by  setting
$$s(h,\bA)=\left(\grad_\bA (h|_{\bC^r}), 
\left. (\hess_\bA  h) \right|_{\nu(\bC^r)}\right).$$ 
The abundance  condition implies that $s$ is a submersion along 
$\{ h_0 \} \times \bC^r$ (see Proposition 3.4  of \cite{bh}).

Hence there is an open set $V \subset \FF$ containing $h_0$ 
and a subset $V_1 \subset V$ of second  category   
such that $h\in V_1$ implies 
$h|_{\bC^r}$ satisfies Definition \ref{equiv morse}. 
 It follows that there is an $SU(3)$ invariant
neighborhood $\tau$ of $\bC^r$ in $\bC$ and an open neighborhood
$V_2$ of $h$ such that $h' \in V_2$ implies  $\crit(h'|_\tau) \subset 
\bC^r.$ 
Consider the compact subset
$C^0 \subset C^*$ obtained by taking the quotient of
$\bC \setminus \tau'$ under $SU(3),$ where  $\tau' \subset \tau$
is some smaller invariant tubular neighborhood.
Repeating the argument given just before Lemma \ref{morsepert}
with $C$ replaced by $C^0$ shows that
there is a second category subset of $V_3 \subset V_2$
such that  $h' \in V_3$ implies that
$h'|_{C^0}$ satisfies   Definition \ref{equiv morse}
as well. This shows that  $h|_\bC$ is equivariantly Morse
for generic $h \in \FF$ near $h_0$.

\begin{lem}  \label{bott morse pert} 
Suppose $\bC \subset \bM_{h_0}$ is a nondegenerate critical submanifold
and $f$ is an admissible function
such that $f|_{\bC}$ is equivariantly Morse. 
Set $h_t=h_{0}+tf$ and let $C \subset \MM_{h_0}$
be the image of $\bC$ under $\bM_{h_0} \to \MM_{h_0}.$
Then there is an open set $ U \subset \BB $ containing
$C$ and an $\ep >0$ such that, for every $0<t<\ep$,
$ \OO_{t} := \MM_{h_t} \cap  U$ is a regular subset of $\MM_{h_t}$.
Let $\bO_t$ be the preimage of $\OO_t$ under $\pi:\bB \to \BB.$
There is bijection $\varphi_t : \crit(f |_C) \to  \OO_t$ which lifts to 
an
$SU(3)$ equivariant diffeomorphism 
$\tilde{\varphi}_{t}: \crit (f|_{\bC}) \to \bO_{t}$. 
Given a smooth family $A_t$ with $ \bbl A_{0} \bbr \in \crit(f|_{\bC})$
and $\bbl A_t \bbr  = \tilde{\varphi}_t(\bbl A_0 \bbr)$, then
for  $0<  t < \ep$,   
\begin{equation} \label{spectral}
{\rm SF} (\th, A_t) = {\rm SF} (\th, A_{0}) + \ind_{ \bbl A_0 \bbr}(f),
\end{equation}
where $\ind_{\bbl A_0 \bbr}(f)$ is the Morse index of the 
critical point $\bbl A_0 \bbr$ of $f|_\bC$.
If, in addition, $A_0$ and $A_t$ are reducible, then (\ref{spectral})  
holds for the $\hh$ and $\hhp$ components separately:
\begin{equation} \label{spectral-b}
\begin{split}
{\rm SF}_{\hh}(\th, A_t) &={\rm SF}_{\hh}(\th, A_{0})+\ind^t_{\bbl A_0 
\bbr}(f),\\ 
{\rm SF}_{\hhp}(\th,A_t) &= {\rm SF}_{\hhp}(\th,A_{0})+\ind^n_{\bbl A_0\bbr}(f),
\end{split}
\end{equation}
where $\ind^t_{\bbl A_0 \bbr} (f )$ and $\ind^n_{ \bbl A_{0} \bbr} (f)$
are the indices of  $\hess_{\bbl A_0 \bbr}  (f|_\bC)$ 
in the directions tangent and normal to $\bC^r$ in $\bC,$ respectively.
\end{lem}

\begin{proof}  Since the argument is nearly identical to the proof of 
Proposition \ref{morsepert}, we only explain the modifications one
needs to make. 
The tangent space to the gauge group $\GG$ at the identity
is given by the space 
of 0-forms completed 
in the $L^2_2$ norm. Therefore, the tangent space to the subgroup $\GG_0 \subset \GG$
of based gauge transformations is the subspace
$$\Om^0_0 =\{ \xi \in L^2_2( \Om^0(M;su(3))) \mid \xi(x_0)=0\}$$ 
consisting of 0-forms vanishing at the basepoint.
Consider the bundle $\bL$ whose fiber above $(A,h)$ is
$$\bL_{A,h}=\{  (\xi,a )\in \JJ_{A,h} \mid \xi=0, a\perp d_A (\Om^0_0) \},$$
and denote again by $\bL$  the induced bundle on the quotient
$\bB\times\FF$.  Notice that the fiber $\bL_{A,h}$ contains $\ker d_A^*$ as a subspace of 
codimension $8-\dim \Gamma_A$.  We regard $\bL_{A,h}$ as the tangent 
space of $\bB$ at $\bbl A \bbr$.  
 
By restricting and projecting $K_{A,h}$, we obtain an operator $\bK $ on 
$\bL.$  This operator agrees with $K_{A,h}$ on $\ker d_A ^*$ and vanishes on
the orthogonal 
complement to $\ker d_A ^*$ in $\bL_{A,h}$, which is 
just the tangent space to the orbit of the 
residual $SU(n)$ action.
Choose open subsets $\bU \subset \bB$ containing $\bC$ and
$V \subset \FF$ containing $h_0$ and
define $\la_0$ as in (\ref{lam}), with $\hK$ replaced by $\bK$.
Over $\bU \times V,$ decompose $\bL = \bL_0 \oplus \bL_1$
into the two eigenbundles as in (\ref{big and small}). 
Let   $\wQ:\bC \times (-\ep,\ep) \to \bL_0$ be
the analog of the map $Q$ from before.

The only substantial difference is that
now $f|_{\bC}$ is not Morse but rather 
equivariantly Morse. This implies 
that $f$ induces Morse functions on $C^*$ and $C^r$
with only finitely many critical points.
The  argument from Proposition
\ref{morsepert} which produced the map $\hQ$ on $C\times (-\ep, \ep)$ 
can also be applied here and results in equivariant maps   
$\bC^{*}\times (-\ep, \ep) \to \bL_0$ and $\bC^{r}\times (-\ep,\ep) \to 
\bL_0$ 
whose zero sets together coincide with that of $\wQ$.  Reducing modulo $SU(n)$,
we obtain 
1-dimensional (product) cobordisms in $\BB^{*}$ and $\BB^{r}$ 
which we follow to define the map $\varphi_t$.  The
preimages of the cobordisms under $\pi:\bB \to \BB$
are equivariant product cobordisms in $\bB$.
Nondegeneracy of $\hess f$ in the normal direction to $\bC^{r}$ 
guarantees that there are no irreducible orbits in $\wQ^{-1}(0)$ nearby,
and the claims about the spectral flow follow as in the previous 
case.
\end{proof}

The following proposition
applies to components of type (ii) and determines their contribution
to $\la_{SU(3)}(X_1 \# X_2).$

\begin{prop}   \label{red comp}
Suppose $h_0$ is a small perturbation 
and $C \subset \MM_{h_0}$ is a connected component satisfying
the following conditions:
\begin{enumerate}
\item[(i)] For each $[A] \in C$,
the isotropy group $\Ga_A$ is isomorphic to either $\ZZ_3$ or $U(1)$.
\item[(ii)] The lift $\bC$ of
$C$ under the projection $\pi:\bB \to \BB$ is a nondegenerate critical 
submanifold.
\item[(iii)] Both $C^*$ and $C^r$ are connected. 
\end{enumerate}
Choose connections $A_0,B_0$ with $[A_0] \in C^*$ and
$[B_0] \in C^r$. 
Then the contribution of 
 $C$ to $\la_{SU(3)}(X)$ is 
\begin{equation} \label{relEuler}
(-1)^{{\rm SF}(\th,A_0)} \chi(C, C^r)-\tfrac{1}{2} (-1)^{{\rm SF}(\th,B_0)} \chi(C^r)
\left( {\rm SF}_\hhp(\th,B_0) - 4  \CS(\widehat{B}_0)+2\right) 
\end{equation}
 where $\widehat{B}_0$ is a  flat, reducible connection close to $B_0$.
\end{prop}

\noindent
{\it Remark. \;}  We do not assume  
$X$ is a connected sum 
in either Proposition \ref{irred comp} or \ref{red comp} as there
may be other interesting applications of these results, e.g.,
to components of the flat moduli space
of positive dimension. 
Condition (iii)  
holds for components  $C$ arising from 
connected sums   but is not an essential hypothesis.
For example, if $C^r$ is not connected, then decompose it
into its connected components
$$C^r = \bigcup_{i=1}^m C^r_i$$ 
and choose $B_i \in \AA$ with
$[B_i] \in C^r_i$ for $i=1,\ldots, m$.
 Then the correct statement is obtained by replacing 
 (\ref{relEuler}) by 
$$(-1)^{{\rm SF}(\th,A_0)} \chi(C, C^r)-\tfrac{1}{2} \sum_{i=1}^m (-1)^{{\rm SF}(\th,B_i)} 
\chi(C^r_i) \left( {\rm SF}_\hhp(\th,B_i) - 4  \CS(\widehat{B}_i)+2\right).$$

\begin{proof}
We first show that (\ref{relEuler}) is independent
of the choices of $A_0,B_0$ and $\widehat{B}_0$.
The argument
of Theorem 5.1 of \cite{bh}  shows that  
(\ref{relEuler}) depends only on  the gauge orbits
$[A_0], [B_0] \in C$ and not on their gauge representatives.
That argument also shows that (\ref{relEuler}) 
is independent of the choice of $\widehat{B}_0$.
So, it suffices to  show that (\ref{relEuler})
is independent of the choice of $[A_0] \in C^*$ and $[B_0] \in C^r$.

The Lie group $SU(3)$ acts smoothly on $\bC$, and  hence Corollary  
VI.2.5 of \cite{bredon} implies $\bC^r$
is a smooth submanifold of $\bC$. Since
$PU(3)= SU(3)/\ZZ_3$ acts freely on $\bC^*,$
the quotient $C^*$ is also smooth.   
Thus the dimension of the kernel of $\hess_A (CS+h_0)$ 
is constant as a function of 
$[A] \in C^*$ (the tangent space of $C^*$ at $[A]$
can be identified with 
the space of zero modes of the Hessian).  
The same is true
of the signature operator 
$$K_A:\Om^{0+1}(X,su(3)) \lto \Om^{0+1}(X,su(3)),$$ 
since it is just
 the Hessian enlarged by putting
$d_A: \Om^0(X,su(3)) \to \Om^1(X,su(3))$ and
its adjoint
$d^* _A : \Om^1(X,su(3)) \to \Om^0(X,su(3))$ in opposite off-diagonal blocks.

Given $[A_0], [A'_0 ] \in C^*$, there is by (iii) a path in $C^*$ 
from $[A_0]$ to $[A'_0]$ which 
we lift to a path $A_t$ of irreducible
connections from $A_0$ to $A_1= g \cdot A'_0$, where $g \in \GG$.
Since none of the eigenvalues
of $K_A$
cross zero along $A_t$, it follows that
  ${\rm SF}(\th,A_0) = {\rm SF}(\th, A_1)$.
This proves 
(\ref{relEuler})  is independent of the choice
of $[A_0] \in C^*$.

To prove (\ref{relEuler})  is independent of the choice 
of $[B_0] \in C^r,$ choose a lift $\bbl B_0 \bbr \in \bC^r$ of
$[B_0]$ and   decompose the  
tangent space of $\bC$ at $\bbl B_0 \bbr$
into the subspaces of
vectors tangent to $\bC^r$ and vectors normal 
to $\bC^r$ in $\bC$. 
Now $C^r$ connected implies $\bC^r$ is connected,
and hence the dimension of the kernel of  
$\hess_B (CS+h_0)$ is constant as a function
of $\bbl B \bbr \in \bC^r$. 
The same is true for the restriction
$$ \hess_B (CS+h_0) \vert_{\Om^1(X;\hhp)}$$
because its kernel can be identified with the normal bundle
of $\bC^r$ in $\bC$. 
Similar statements hold for  the signature operator $K_B$
and its restriction $K_B \vert_{\Om^{0+1}(X;\hhp)}$
(notice that $H^0_B(X;su(3)) = \RR$ and $H^0_B(X; \hhp) = 0$
for $\bbl B \bbr \in \bC^r$).

Given $[B_0], [B'_0] \in C^r,$ there is
a path in $C^r$  from $[B_0]$ to $[B'_0]$
which we lift to a path $B_t$ of reducible
connections from $B_0$ to $B_1=g\cdot B'_0$. 
Since none of the eigenvalues
of $K_B$ or its restriction $K_B \vert_{\Om^{0+1}(X;\hhp)}$
cross zero along $B_t$, it follows that
  ${\rm SF}(\th,B_0) = {\rm SF}(\th, B_1)$
and ${\rm SF}_\hhp(\th,B_0) = {\rm SF}_\hhp(\th,B_1)$.
This proves  
that (\ref{relEuler})  is independent of the choice
of $[B] \in C^r$.

To compute the contribution of $C$ to $\la_{SU(3)}(X)$,
we choose an admissible function $f$
so that $f|_\bC$ is equivariantly Morse and consider the
parameterized moduli space 
$$W = \bigcup_{0\leq t \leq t_0} \MM_{h_t}\times \{t\}$$ for the 1-parameter
family of perturbations $h_t=h_0 + tf.$  For $t_0$ small, 
$W$ is a union of 
connected components corresponding to the connected 
components of $h_0$.  Let $U$ be the component of $W$ containing 
$C\times \{ 0 \}$, and let $U_t$ denote the ``t-slice'' $ U\cap ( \MM_{h_t}\times \{t\} ) $.
  
Then, by definition, the contribution of $C$ to $\la_{SU(3)}(X)$ 
is the sum
\begin{equation} \label{izzy}
\sum_{[A] \in U^*_{t_0}} (-1)^{{\rm SF}(\th,A)}
- \tfrac{1}{2} \sum_{[B] \in U^r_{t_0}} 
(-1)^{{\rm SF}(\th,B)}({\rm SF}_\hhp(\th,B) - 4 \CS(\widehat{B})+2)
\end{equation}
where $t_0$ is a small positive number  
and $U_t = U^*_t \cup U^r_t$ is the decomposition into
irreducible and reducible gauge orbits.

From equation (\ref{spectral}) of Lemma \ref{bott morse pert}, 
it follows that
\begin{equation} \label{izzypop}
\begin{split}
\sum_{[A] \in U^*_{t_0}} (-1)^{{\rm SF}(\th,A)}
&= \sum_{[A] \in \crit(f|_{C^*})} (-1)^{{\rm SF}(\th,A)}(-1)^{\ind_{\bbl A\bbr }(f)} \\
&=(-1)^{SF(\th, A_0)} \sum_{[A]\in \crit(f|_{C^*}) } (-1)^{\ind_{[A]}(f)}.
\end{split}
\end{equation}
This uses the previously
established fact that $(-1)^{{\rm SF}(\th,A)}= (-1)^{{\rm SF}(\th,A_0)}$ 
for all $[A] \in C^*$, together with the observation that
the Morse index of $f$ at $\bbl A\bbr \in \bC^*$ equals that 
of the induced function $f$ on $C^*$ at $[A]$.   

Similarly, from equation (\ref{spectral-b}) of Lemma \ref{bott morse pert},
 it follows that
 \begin{align}
 \sum_{[B] \in U^r_{t_0}} &
(-1)^{{\rm SF}(\th,B)}({\rm SF}_\hhp(\th,B) - 4 \CS(\widehat{B})+2)  \nonumber \\
& = \sum_{[B] \in \crit(f|_{C^r})} (-1)^{{\rm SF}(\th,B)}
(-1)^{\ind_{\bbl B \bbr}(f)} (
  \ind^n_{\bbl B \bbr}(f)+{\rm SF}_\hhp(\th,B) 
- 4 \CS(\widehat{B})+2)  \nonumber \\
&=  \sum_{[B] \in \crit(f|_{C^r})} (-1)^{{\rm SF}(\th,B_0)}
(-1)^{\ind^t_{\bbl B \bbr}(f)} (\ind^n_{\bbl B \bbr}(f) 
+ {\rm SF}_\hhp(\th,B_0)- 4 \CS(\widehat{B}_0)+2)  \nonumber \\
\begin{split}
&=  (-1)^{{\rm SF}(\th,B_0)} ({\rm SF}_\hhp(\th,B_0)- 4 \CS(\widehat{B}_0)+2)
\sum_{[B] \in \crit(f|_{C^r})}
(-1)^{\ind^t_{\bbl B \bbr}(f)}   \\
& \qquad -  (-1)^{{\rm SF}(\th,A_0)} \sum_{[B] \in \crit(f|_{C^r})} 
(-1)^{\ind^t_{\bbl B \bbr}(f)} ( \ind^n_{\bbl B \bbr}(f)). \label{izzypup}
\end{split}
\end{align}
The second step follows since $(-1)^{SF(\th,B)}$ and
$ {\rm SF}_\hhp(\th,B)- 4 \CS(\widehat{B})$  are independent of
$[B] \in C^r$ and since   
$\ind^n_{\bbl B \bbr}(f)$ is even. The last step
is justified by the following  lemma.

\begin{lem} For all $[A] \in C^*$ and all $[B] \in C^r$, 
 $(-1)^{{\rm SF}(\th,B)} = (-1)^{{\rm SF}(\th,A) +1}$. \label{platonic relationship}
\end{lem}
\begin{proof}
To prove the lemma, suppose
$\be_t$ is a 1-parameter family in $\AA$ with $[\be_0] \in C^r$  and
$[\be_t] \in C^*$   for $t>0$. Then
\begin{align*}
\dim H^0_{\be_1}(X;su(3)) &= \dim H^0_{\be_0}(X;su(3))-1, \hbox{ and} \\
\dim H^1_{\be_1,h_0}(X;su(3)) &= \dim H^1_{\be_0,h_0}(X;su(3))-1.
\end{align*}
Indeed, as $t$ increases from $t=0,$ a pair of eigenvalues of $K_{\be_t,h}$
of equal magnitude and opposite sign leave zero.
This proves that
${\rm SF}(\th,\be_0) = {\rm SF}(\th,\be_1) -1$. It also proves the claim
since, as we have already seen,  $ (-1)^{{\rm SF}(\th,B)}$ is independent
of $[B] \in C^r$ and $ (-1)^{{\rm SF}(\th,A)}$ is independent
of $[A] \in C^*.$\end{proof}

We now complete the proof of Proposition \ref{red comp}.
Substituting equations (\ref{izzypop}) and (\ref{izzypup}) into
(\ref{izzy}), we see that the contribution of $C$ to $\la_{SU(3)}(X)$ is given
by 
\begin{equation}  \label{first step}
\begin{split}
	& (-1)^{SF(\theta, A_0)} 
		\left[ \sum_{[A] \in \crit(f|_{C^*})} (-1)^{\ind_{[A]} (f)} 
			+ \tfrac{1}{2} \sum_{[B] \in \crit (f|_{C^r})} 
			(-1)^{\ind^t _{\bbl B\bbr} (f)} 
			(\ind^n_{\bbl B\bbr} (f)) \right]
\\ 
&	\qquad 	
-\tfrac{1}{2} \left( SF_{\hh ^\perp} 
(\theta ,B_0) -4 \CS (\widehat{B}_0) + 2 \right) 
(-1)^{SF(\theta, B_0)} \sum_{[B]\in \crit(f|_{C^r})} (-1)^{\ind^t _{\bbl B \bbr} (f)} \\
\end{split}
\end{equation}

Notice that quantity in brackets on the first line of (\ref{first step})  is independent 
of the equivariantly Morse function $f$ on $\bC$. (This follows from 
an argument similar to but simpler than that given in \cite{bh} to 
show that $\lambda_{SU(3)}$ is independent of perturbation.) 
Hence we can compute it using any equivariantly Morse
function we want. Choosing  a function
whose Hessian   in the normal directions to 
$\bC^r$ is positive definite and whose critical values along 
$\bC^*$ are all larger than the values along $\bC^r$, we
see that the quantity in brackets on the first line of (\ref{first step}) equals
the relative Euler characteristic $\chi(C,C^r).$
A standard argument shows that
$$ \sum_{[B]\in \crit(f|_{C^r})} (-1)^{\ind^t _{\bbl B \bbr} (f)}  = \chi(C^r).$$
This proves (\ref{first step}) equals (\ref{relEuler}) and we are done.
\end{proof}

We can now complete the proof of Theorem \ref{connect}.
As explained earlier, the point components in $\MM_{h_0}$
give rise to the first two terms on the right of formula (\ref{c form}).
Further, if $C$ is a connected component of $\MM_{h_0}$ of type (i),
then it contributes algebraically zero to $\la_{SU(3)}(X)$.
This follows from Proposition \ref{irred comp} since
$\chi(C)=0 $ for such $C$. (See Proposition \ref{component types}.
In the case $C \cong N$, this follows simply because $N$ is an 
orientable manifold of odd dimension.)

It remains to determine the contribution 
to $\la_{SU(3)}(X_1 \# X_2)$ of components $C$ of  type (ii).
Our first step will be to  calculate the relative Euler characteristic
$\chi(C,C^r)$. By the exactness property of singular homology,
$$\chi(C,C^r) = \chi(C) - \chi(C^r) = \chi(C),$$
where the last step follows from the fact that
$C^r \cong SO(3)$, which is well-known in $SU(2)$
gauge theory.  (See~p.134, \cite{wli}.)
Our computation of $\chi(C)$ utilizes the following 
description of $C$ as the quotient of a certain $U(1)$ action on $N.$

Recall that $N= SU(3)/U(1)$ is our model for
fibers of $\bB \to \BB$ above reducible orbits.
In terms of a reducible $SU(3)$
representation $\varrho$ of $\pi_1(X)$, $N$ is just
the adjoint orbit of $\varrho$, namely
points in $N$ correspond to $SU(3)$ representations conjugate to 
$\varrho.$
Because these representations are all reducible, associated to each point
in $N$ there is a canonical 1-dimensional subspace of $\CC^3$ given by
the invariant linear subspace of the corresponding representation.
This defines a map $N \to \CP^2$ which is, in fact a  fibration.
The fiber above $[0,0,1] \in \CP^2$ consists of $SU(3)$ 
representations $\vartheta$ conjugate to $\varrho$
with $\im(\vartheta) \subset SU(2) \times 1.$
The two irreducible $SU(2)$ representations $\varrho'$ and $\vartheta'$
associated to $\varrho$ and $\vartheta$
are conjugate, and hence the fiber  of $N \to \CP^2$ is
$SO(3),$ the adjoint $SU(2)$ orbit of $\varrho'.$
 
In general, define  
$$\Ga_\varrho= \{ g \in SU(3) \mid g \varrho g^{-1} = \varrho \}$$
and recall that $\varrho$ is reducible and nontrivial if and only if 
$\Ga_\varrho \cong U(1).$
Suppose  $\varrho_1$ and $\varrho_2$ are 
nontrivial reducible $SU(3)$ representations
of $X_1$ and $X_2,$ respectively.
Then $C$ consists of the conjugacy classes of representations
$\varrho$ of $X_1 \# X_2$ such that the restriction of
$\varrho$ to $\pi_1 (X_i)$ is conjugate to $\varrho_i$ for $i=1,2.$
Proposition \ref{component types} shows that   
$$\bC = SU(3)/\Ga_{\varrho_1} \times SU(3)/\Ga_{\varrho_2} \cong N \times 
N,$$
and  by fixing the first factor,
it follows that $C$ is the quotient of the second factor by
the induced action of $\Ga_{\varrho_1} \cong U(1).$
If $\varrho_1$ is chosen with image contained in $SU(2) \times 1$, 
then $\Ga_{\varrho_1}$ is
simply the $U(1)$ subgroup described in (\ref{U1 subgroup}).
The 
subgroup of this group consisting of cube roots of 1 acts trivially.
 
The $U(1)$ action descends to the base of the fibration $\pi:N \to 
\CP^2$,
where it acts by $[x,y,z] \mapsto [ux,uy,u^{-2} z]$ and
has fixed point set $\{ [0,0,1]\} \cup \{[x,y,0] \}= \{pt\} \cup \CP^1$.
Notice that  $\pi^{-1}([0,0,1]) = C^r$ and set  
$B_1=\CP^2 \setminus \{[0,0,1]\}$ and $B_2 = \CP^2 \setminus \CP^1.$
Define $C_i=\pi^{-1}(B_i)/U(1)$ and observe that
$C=C_1 \cup C_2$ and $C^r \subset C_2.$
The  Mayer-Vietoris sequence gives that
$$\chi(C) = \chi(C_1)+\chi(C_2)-\chi(C_1 \cap C_2).$$
However, $U(1)/\ZZ_3$ acts freely on  
$B_2 \setminus \{[0,0,1]\} \cong \CC^2 \setminus \{0\}$
and trivially on the fiber above $[0,0,1]$,
and hence $C_2$ is an $SO(3)$ bundle over $B_2/U(1).$
Thus, $\chi(C_2)=0.$ Similarly, $\chi(C_1 \cap C_2) =0.$
 
Now $B_1$ certainly retracts to $\CP^1,$ and
we claim that $C_1$ also retracts to $C_0=\pi^{-1}(\CP^1)/U(1)$.
This follows by considering the  
$U(1)$ action on the fibers above $[x,y,0] \in \CP^1.$ 
For example, take $p=[1,0,0]$, the north pole. 
If $\varrho_2 \in \pi^{-1}(p),$ then $\im(\varrho_2) \subset H$ where
$$H =  \left\{ \left. \left( 
\begin{array}{crr}1&0&0\\
0&a&b\\
0&-\bar{b} &\bar{a} \end{array}
 \right)
\right| a \bar{a} + b \bar{b}=1. \right\}.$$
In this case,  $U(1)$ acts on $\pi^{-1}(p) = SO(3)$ in the standard
way by rotation of the off-diagonal entries and has quotient $ S^2.$
Hence $\pi^{-1}(p) / U(1) \cong S^2$ and nearby, the $SO(3)$
fibers in $C^*$ retract to the $S^2$ fibers in $\pi^{-1}(\CP^1)/U(1)$
via the cone structure.

The same is true for $[x,y,0] \in \CP^1.$ 
For suppose $x$ and $y$ are complex numbers satisfying
$x \bar{x} + y \bar{y} =1$ and suppose that
$$\al = \left( 
\begin{array}{rrr} x&  -\bar{y}  &0\\
 y &  \bar{x} & 0\\
 0 & 0 & 1\end{array} \right).$$
Then $\al(p) = [x,y,0]$ and $\im(\varrho_2) \subset \al H \al^{-1}$ 
whenever
$\varrho_2$ is a reducible $SU(3)$ representation with  
invariant linear subspace $[x,y,0].$ 
Since the $U(1)$ action commutes with multiplication by $\al,$ 
it acts on  
$\al H \al^{-1}$ in the same way as it did on
$H$. Thus $\pi^{-1}([x,y,0])/ U(1) \cong S^2$ and
there is a fibration $ C_0 \to \CP^1$ with fiber  $S^2$.
Hence  $\chi(C_0) = 4$
  and we conclude that $\chi(C, C^r) = 4.$

Since $\chi(C^r)=0$ for components of type (ii), 
all the terms involving $B$ in Proposition \ref{red comp}
vanish and it follows that each such component
contributes $(-1)^{{\rm SF}(\th,A)}  \chi(C,C^r)$ to 
$\la_{SU(3)}(X_1 \# X_2)$,
where $A = A_1 \#_\si A_2$ is chosen so that
$[A] \in C^*.$ In order to compute $ {\rm SF}(\th,A) \mod 2,$
it is convenient to set $B=A_1 \#_\tau A_2$ with
$[B] \in C^r$. Since $B$ is reducible,  
$\dim H^0_B(X;su(3)) = 1$ and we compute that
$\dim H^1_B(X;su(3)) = 7$
(this uses the splitting $ su(3) = 
su(2) \oplus \CC^2 \oplus \RR$ together with the facts:
$\dim H^1_B(X;su(2)) = \dim C^r = 3$,
$H^1_B(X;\RR)=0$, and
$ H^1_B(X;\CC^2)=4$).
On the other hand, if $[A] \in C^*,$
then $ H^0_A(X;su(3)) = 0$ and 
$\dim H^1_A(X;su(3)) = \dim C^* = 6.$
By Lemma \ref{platonic relationship}, 

\begin{equation} \label{stepone}
{\rm SF}(\th,A)  \equiv {\rm SF}(\th,B)-1 \mod 2.
\end{equation}

Using the splitting $su(3)= su(2) \oplus \CC^2 \oplus \RR$,
and applying
 equation (\ref{splitit}) to the $su(2)$ component of  ${\rm SF}(\th, B)$,
we see that
\begin{eqnarray}
{\rm SF}(\th, B) &=& {\rm SF}(\th, A_1 \#_\tau A_2) \nonumber \\
&=& {\rm SF}_{su(2)}(\th, A_1 \#_\tau A_2) + {\rm SF}_{\CC^2}(\th, A_1 \#_\tau A_2) +
 {\rm SF}_{\RR}(\th, A_1 \#_\tau A_2) \nonumber \\
 &\equiv& {\rm SF}_{su(2)}(\th_1, A_1) +  {\rm SF}_{su(2)}(\th_2,A_2) -1 \mod 2 \nonumber \\
 &\equiv& {\rm SF} (\th_1, A_1) +  {\rm SF}_{su(2)}(\th_2,A_2) - 1 \mod 2. \label{steptwo}
 \end{eqnarray}
 The third and fourth steps
 follow because all the $\CC^2$ spectral flows are even and all the
 $\RR$ spectral flows equal $-1$ (the coefficients are untwisted,
 $X$ is a homology sphere, and we use the
 $(-\epsilon, \epsilon)$ convention for computing spectral flows).
Combining equations (\ref{stepone}) and (\ref{steptwo}), we conclude that
$${\rm SF}( \th, A) \equiv {\rm SF}(\th,B)-1 \equiv {\rm SF}(\th_1,A_1)+ {\rm SF}(\th_2,A_2) \mod 2.$$

This, together with the above computation  of $\chi(C,C^r)$,
implies 
\begin{eqnarray*}
\sum (-1)^{{\rm SF}(\th,A)}  \chi(C,C^r)
&=& 4\!\!\!\!\! {\displaystyle \sum_{[A_1] \in \MM^r_{h_1}(X_1)}} 
 \!\!\!\!\! (-1)^{{\rm SF}(\th_1, A_1)}
{\displaystyle \sum_{[A_2] \in \MM^r_{h_2}(X_2)}} 
\!\!\!\!\!\! (-1)^{{\rm SF}(\th_2, A_2)} \\
&= & 4 \, \la_{SU(2)}(X_1) \; \la_{SU(2)}(X_2),
\end{eqnarray*}
where the first sum is over all components $C \subset \MM_{h_0}$
of type (ii) and $[A] \in C^*$.
Recall that $h_0 = h_1 + h_2$ and
$\MM_{h_i}(X_i)$ is regular for $i=1,2$.
 This completes the proof of Theorem \ref{connect}.
\bigskip

{\bf Acknowledgements.}
{The first named author was partially supported by 
a Seed Grant from Ohio State University and the second named 
author was partially supported by a Research Grant from Swarthmore 
College.
Both authors wish to thank Thomas Hunter, Paul Kirk and Eric 
Klassen
for helpful discussions.



\end{document}